\documentclass[aos,MSNbibl,nameyear,dvips]{arximspdf}

\doi{10.1214/12-AOS972} 
\volume{40}
\issue{1}
\pubyear{2012}
\firstpage{593}
\lastpage{608}

\makeatletter

\newcommand{\bigtimes}{\mathop{\mbox{\fontsize{17}{17}\selectfont{$\!\times$}}}}

\newcommand{\cal}{\mathcal}

\renewcommand{\mid}{\,|}

\newcommand{\clos}{\operatorname{cl}}
\newcommand{\bound}{\operatorname{bd}}
\newcommand{\defeq}{:=}
\newcommand{\transp}{^{\mathrm T}}
\newcommand{\reals}{\mathbb{R}}
\newcommand{\E}{{\mathrm{E}}}

\newcommand{\bq}{Q}
\newcommand{\bqq}{\mathbf{q}}
\newcommand{\bp}{P}
\newcommand{\bpp}{\mathbf{p}}
\newcommand{\bS}{\mathbf{S}}

\newcommand{\balpha}{{\bolds\alpha}}
\newcommand{\bbeta}{{\bolds\beta}}
\newcommand{\cX}{{\cal X}}
\newcommand{\cA}{{\cal A}}
\newcommand{\cP}{{\cal P}}

\newcommand{\neighbours}[1]{\mathrm{ne}(#1)}

\newcommand{\cN}{\mathcal{N}}
\newcommand{\cB}{\mathcal{B}}
\newcommand{\graf}{{\mathcal{G}}}

\newtheorem{theorem}{Theorem}[section]
\newtheorem{lemma}[theorem]{Lemma}
\newtheorem{cor}[theorem]{Corollary}

\newproclaim{define}{Definition}[section]
\newproclaim{cond}{Condition}[section]
\newproclaim{expl}{Example}[section]

\makeatother

\begin{document}
\begin{frontmatter}

\title{Proper local scoring rules on discrete sample~spaces}
\runtitle{Discrete local scoring}

\begin{aug}
\author[A]{\fnms{A. Philip} \snm{Dawid}\ead[label=e1]{A.P.Dawid@statslab.cam.ac.uk}},
\author[B]{\fnms{Steffen} \snm{Lauritzen}\corref{}\ead[label=e2]{steffen@stats.ox.ac.uk}}
\and
\author[C]{\fnms{Matthew} \snm{Parry}\thanksref{t1}\ead[label=e3]{mparry@maths.otago.ac.nz}}
\runauthor{A. P. Dawid, S. Lauritzen and M. Parry}
\affiliation{University of Cambridge, University of Oxford
and University of Otago}
\address[A]{A. P. Dawid\\
Statistical Laboratory\\
Centre for Mathematical Sciences\\
University of Cambridge\\
Wilberforce Road\\
Cambridge CB3 0WB\\
United Kingdom\\
\printead{e1}}
\address[B]{S. Lauritzen\\
Department of Statistics \\
University of Oxford\\
South Parks Road\\
Oxford OX1 3TG\\
United Kingdom\\
\printead{e2}}
\address[C]{M. Parry\\
Department of Mathematics\\
\quad and Statistics\\
University of Otago\\
P.O. Box 56\\
Dunedin 9054\\
New Zealand\\
\printead{e3}} 
\end{aug}

\thankstext{t1}{Supported by EPSRC Statistics Mobility Fellowship
EP/E009670.}

\received{\smonth{4} \syear{2011}}
\revised{\smonth{1} \syear{2012}}

%
\begin{abstract}
A~\textit{scoring rule} is a~loss function measuring the
quality of a~quoted probability distribution $\bq$ for a~random
variable $X$, in the light of the realized outcome $x$ of $X$; it is
\textit{proper} if the expected score, under any distribution $\bp$
for $X$, is minimized by quoting $\bq= \bp$. Using the fact that
any differentiable
proper scoring rule on a~finite sample space ${\cal X}$ is the
gradient of a~concave homogeneous
function, we consider when such a~rule can be
\textit{local} in the sense of depending only on the probabilities
quoted for points in a~nominated neighborhood of $x$. Under mild
conditions, we characterize such a~proper local scoring rule in
terms of a~collection of homogeneous functions on the cliques of an
undirected graph on the space $\cX$. A~useful property of such rules
is that the
quoted distribution $\bq$ need only be known up to a~scale factor.
Examples of the use of such scoring rules include Besag's
pseudo-likelihood and Hyv\"arinen's
method of ratio matching.
\end{abstract}

%

\end{frontmatter}

\section{Introduction}
\label{secintro}

Let $\mathcal{X}$ be a~finite set, let $\cA$ be the set of real
vectors $\balpha= (\alpha_x\dvtx x \in\cX)$ with each $\alpha_x > 0$,
and let $\mathcal{P} = \{\bpp\in\cA\dvtx \sum_x p_x = 1\}$ be the set of
such vectors corresponding to strictly positive probability
distributions on~$\mathcal{X}$. We will use $\bp$ for the
distribution determined by $\bpp$ (similarly $Q$ for~$\bqq$), and
generally do not distinguish between them. For $\balpha\in\cA$, $C
\subseteq\cX$ we write $\balpha_C \defeq(\alpha_x \dvtx x \in C)$, and
similarly $\bpp_C$.

Consider a~game between Forecaster and Nature, where Forecaster quotes
a~distribution $\bq\in\cP$ as representing his uncertainty about a~quantity $X$ taking values in $\cX$, and Nature then reveals $X=x$. A~\textit{scoring rule} [see, e.g., \citet{apdencprobfore}] is a~function $S\dvtx \cX\times\cP\rightarrow\reals$. The interpretation
is that~$S(x, \bq)$ measures the loss suffered by Forecaster for the
above outcome of the game.

For $\bp\in\cP$ we define $S(\bp,\bq) \defeq\sum_x p_x S(x, \bq)$,
the expected score when Forecaster quotes $\bq$, and Nature generates
$X$ from $\bp$. The scoring rule $S$ is \textit{proper} if always
$S(\bp,\bq) \geq S(\bp,\bp)$, so that it is always optimal to quote a~distribution $\bq$ matching the real uncertainty $\bp$; $S$ is
\textit{strictly proper} if furthermore $S(\bp,\bq) > S(\bp,\bp)$ when
$\bq\neq\bp$.

The generalized \textit{entropy function}, or \textit{uncertainty
function}, $H\dvtx \cP\rightarrow\reals$, associated with a~proper
scoring rule $S$ is given by $H(\bp) \defeq S(\bp,\bp)$. Then $H$ is
a~concave function on $\cP$. We also introduce the associated
\textit{divergence} or \textit{discrepancy} function $d\dvtx \cP\times\cP
\rightarrow\reals$, where $d(\bp,\bq):=S(\bp,\bq)-H(\bp)$. Then
$d(\bp, \bq) \geq0$, with equality if $\bq= \bp$ (and only in this
case if $S$ is strictly proper).

As well as being of intrinsic interest, proper scoring rules have a~range of applications. For example, if ${\cal Q} = \{\bq_\theta:
\theta\in\Theta\}$ is a~smooth parametric statistical model, we might
estimate $\theta$, based on a~random sample $(x_1, \ldots, x_n)$, by
minimizing the \textit{empirical discrepancy},
$d(\hat\bp_n,\bq_\theta)$, where $\hat\bp_n$ is the empirical
distribution of the data. This is equivalent to minimizing
$\sum_{i=1}^n S(x_i, \bq_\theta)$. Implementing this by setting the
derivative of this criterion to $0$ will yield an unbiased estimating
equation [\citet{Dawid2005ge}, \citet{apd07aism}], from which we, under
suitable smoothness assumptions, can deduce
statistical properties of the associated
estimator such as consistency and asymptotic normality. For the well-known
\textit{logarithmic score} $S(x,Q)=-\log q_x$ this procedure leads to
the maximum likelihood
estimator, but it is of interest to use other scoring rules and
estimators, for example,
because they can lead to greatly simplified calculations. We illustrate this
in Section~\ref{secex} for capture--recapture experiments and pseudo-likelihood
estimation for image analysis; see also \citet{czadoetal09} for a~range of other applications of proper scoring rules to discrete data.

\citet{mfpapdsllplsr} investigated when a~proper scoring rule with
$\cX$ an interval on the real line can be \textit{local} in the sense
that $S(x,Q)$ depends on $Q$ only through the value $f(x)$ of the
density $f$ of $Q$ and the values $f^{(k)}(x)$ of a~finite number of
the derivatives of $f$ at the realized outcome $x$ and
\citet{ehmgneitingorder2} studied rules with $k=2$ in further detail.
It was shown in \citet{mfpapdsllplsr} that any proper local scoring
rule is a~linear combination of the logarithmic score and what was
termed a~\textit{key local scoring rule}; and that any such key local
scoring rule is \textit{$0$-homogeneous} in the sense that $S(x,Q)$ can
be evaluated when $f$ is only known up to proportionality. The results
in this article for discrete sample spaces parallel these. However, in
the case of discrete $\cX$ we have to redefine locality using a~neighborhood structure on the space $\cX$, and use somewhat different
techniques of proof.

The organization of the paper is as follows. In Section~\ref{secconstruct} we
review results from \citet{hendricksonbuehler71} characterizing
proper scoring rules as supergradients of concave functions.

In Section~\ref{seclocalscore} we formally define what it means for a~scoring
rule to be local with respect to a~neighborhood system and show that if
the homogeneous extension of\vadjust{\goodbreak} the scoring rule is local, the
neighborhood system must be determined by an undirected graph. We also
describe a~general additive construction of local scoring rules.
Section~\ref{secex} gives examples of local scoring rules and their use. In
Section~\ref{secconverse} we proceed in parallel to
\citet{mfpapdsllplsr}
by a~variational argument to characterize local scoring rules as
solutions to a~key differential equation and, under an additional
condition on the neighborhood system, we show that such local scoring
rules can be expanded in additive terms indexed by complete subsets of
an undirected graph.\looseness=-1

\section{Homogeneous proper scoring rules}
\label{secconstruct}
Further analysis is facilitated by recasting the problem in terms of
homogeneous functions and using the fundamental characterization of
proper scoring rules given by \citet{mccarthy56} and
\citet{hendricksonbuehler71}.

\subsection{Homogeneous functions}
\label{sechomfun}
A~function $f\dvtx\cA\rightarrow\reals$ is called (\textit{positive})
\textit{homogeneous of order $h$}, or \textit{$h$-homogeneous}, if
%
\begin{equation}
\label{eqhhom}
f(\lambda\balpha) \equiv\lambda^h f(\balpha) \qquad\mbox{for all }
\lambda> 0.
\end{equation}
In this paper we shall only need homogeneity of orders $0$ and $1$.
If $f$ is differentiable,~(\ref{eqhhom}) will hold if and only if $f$
satisfies \textit{Euler's equation}:
%
\begin{equation}
\label{eqeuler}
\sum_x \alpha_x \,\frac{\partial f}{\partial\alpha_x} = h f.
\end{equation}

Even when $f$ is not differentiable, in some circumstances we can
reinterpret~(\ref{eqeuler}) so as to continue to apply.
%
\begin{define}
A~vector $\nabla f(\balpha) \in{\cal A}$ is a~\textit{supergradient}
to $f$ at $\balpha$ if, for all $\bbeta\in{\cal A}$,
\[
f(\balpha) + (\bbeta-\balpha)\transp\nabla f(\balpha) \geq
f(\bbeta).
\]
\end{define}

When $f$ is differentiable at $\balpha$ and has a~supergradient
$\nabla f(\balpha)$ there, it must coincide with the gradient vector
$(\partial f/\partial\alpha_x\dvtx x \in{\cal X})$. Lemma~\ref{lemeuler2}
below and Corollary~\ref{cor1hom} extend Euler's equation~(\ref{eqeuler}) to
homogeneous functions with a~supergradient and are equivalent to
Theorem 2.1 of \citet{hendricksonbuehler71} and subsequent remarks,
so we omit the proofs here.
%
\begin{lemma}
\label{lemeuler2}
Suppose $f$ is $h$-homogeneous, and has a~supergradient $\nabla
f(\balpha)$ at $\balpha$. Then
%
\begin{equation}
\label{eqeuler2}
\balpha\transp\nabla f(\balpha) = h f(\balpha).
\end{equation}
\end{lemma}
%
\begin{cor}
\label{cor1hom}
Suppose $f$ is $1$-homogeneous. Then $\bS$ is a~supergradient of
$f$ at $\balpha$ if and only if
\[
\bbeta\transp\bS\geq f(\bbeta)
\]
for all $\bbeta\in{\cal A}$, with equality when $\bbeta=
\balpha$.\vadjust{\goodbreak}
\end{cor}

By the supporting hyperplane theorem, a~function $f$ is concave on
${\cal A}$ if and only if it has a~supergradient at each
$\balpha\in{\cal A}$ (not necessarily unique if $f$ is not
differentiable at $\balpha$). A~\textit{supergradient function}
$\nabla f$ associates a~specific choice of supergradient $\nabla
f(\balpha)$ with each point $\balpha\in{\cal A}$. If $f$ is
$h$-homogeneous,~(\ref{eqeuler2}) holds at each
$\balpha\in{\cal A}$ for any choice of
supergradient function $\nabla f$.

\subsection{Homogeneous scoring rules}
\label{sechomscore}
Clearly, any scoring rule $S(x, \bp)$ can readily be extended to $\cA
$ by defining $S(x,\balpha) \defeq
S(x, \balpha/\alpha_+)$, where $\alpha_+ \defeq\sum_{y\in\cX}
\alpha_y$. So
extended, $S(x,\balpha)$ is a~$0$-homogeneous function of $\balpha$
for every~$x$ and we say that $S(x,\balpha)$ is a~\textit{$0$-homogeneous scoring rule}.

\citet{mccarthy56} states that a~$0$-homogeneous scoring rule
$S$ is proper if and
only if it can be expressed as the supergradient of a~concave
$1$-homogeneous function $H\dvtx\cA\rightarrow\reals$. This is formally
proved in \citet{hendricksonbuehler71} and stated below in
Theorems
\ref{thmhomprop} and~\ref{thmprophom}.
%
\begin{theorem}
\label{thmhomprop}
$\!\!\!$Suppose $H\dvtx\cA\,{\rightarrow}\,\reals$ is concave and $1$-homogeneous.
Let~$\nabla H$ be a~supergradient of $H$, and, for $x\in{\cal X}$,
$\bpp\in{\cal P}$, define $S(x,\bpp)$ to be the $x$-component of the
vector $\bS(\bpp) \defeq\nabla H(\bpp)$. Then $S$ is a~proper
scoring rule, and the associated entropy at $\bpp$ is $H(\bpp)$.
\end{theorem}

We note that the definition $\bS(\balpha) \defeq\nabla H(\balpha)$
can be used to extend the domain of $S$ from ${\cal X} \times{\cal
P}$ to ${\cal X} \times{\cal A}$. The supergradient function
$\nabla H$ can be taken to be $0$-homogeneous
and then $S(x, \balpha)$ is a~$0$-homogeneous
function of $\balpha$.

For the converse direction, starting with a~scoring rule $S$ defined
on ${\cal X} \times{\cal P}$, we let $S(x,\balpha)$ denote its
$0$-homogeneous extension as described above and let
$\bS(\balpha)$ be the vector with $x$-component $S(x,\balpha)$.
%
\begin{theorem}
\label{thmprophom}
Suppose that $S(x,\balpha)$ is a~$0$-homogeneous proper scoring
rule. Define $H(\balpha) \defeq\balpha\transp\bS(\balpha)$.
Then $H$ is 1-homogeneous and concave, and~$\bS(\balpha)$ is a~supergradient of $H$ at $\balpha$.
\end{theorem}

As a~consequence we obtain the following symmetry relation
for the partial derivatives of any $0$-homogeneous proper scoring rule.
%
\begin{cor}
\label{cord2}
$\!\!\!$If $S$ is a~$0$-homogeneous proper scoring rule, and $S(x,\balpha)$
is continuously differentiable on ${\cal A}$ for each $x\in{\cal
X}$, then
%
\begin{equation}
\label{eqd2}
\frac{\partial S(x, \balpha)}{\partial\alpha_y} =
\frac{\partial S(y, \balpha)}{\partial\alpha_x}.
\end{equation}
\end{cor}
\begin{pf}
In this case $H(\balpha) = \balpha\transp\bS$ is differentiable on
${\cal A}$, so $\bS(\balpha)$ is its gradient. It immediately
follows that $H$ is twice continuously differentiable. Then
(\ref{eqd2}) follows from $\partial^2 H/\partial\alpha_y \,\partial
\alpha_x = \partial^2 H/\partial\alpha_x \,\partial\alpha_y$.
\end{pf}
%
\begin{expl}
Examples of proper scoring rules are the \textit{Brier score} $S(x$, $\bpp)
= \|\bpp\|^2-2p_x$ [\citet{brier50}], where $\|\bpp\|^2= \sum_xp_x^2$,
with $1$-homogeneous entropy function
$H(\balpha)=-\|\balpha\|^2/\alpha_+$ and the \textit{spherical score}
$S(x,\bpp) = -p_x/\|\bpp\|$ with $1$-homogeneous entropy function
$H(\balpha)=-\|\balpha\|$ [\citet{good71}, \citet{apd07aism}].
\end{expl}

We shall say that the entropy function $H$ is \textit{regular} if it is
continuous on~$\cA$ and its closure $\clos{H}$ as a~concave function
[\citet{Rockafellar70}, page~52] is finite on the closed cone
$\bar{\cA}=\{\balpha\dvtx\alpha_x\geq0, x\in\cX\}$. In other words, $H$
is regular if it can be extended by continuity to have finite values
for all $\balpha\in\bar\cA$.

Clearly, since $H(\balpha) = \sum_x \alpha_x S(x, \balpha)$,
$H(\balpha)$ is certainly regular if $S(x, \bp)$ is bounded in
$\bp$ for each $x$.
Both the Brier score and the spherical score satisfy this requirement
and have regular entropy functions, but in general boundedness is not
necessary for regularity.

\section{Local scoring rules}
\label{seclocalscore}
In general, as for the Brier and spherical score, $S(x, \bp)$ will
depend on every element of $\bp$. We are
interested in cases where this is not so.

\subsection{Locality}
\label{seclocal}
Suppose we specify, for each $x\in\cX$, a~set $N_x
\subseteq\cX$ (the \textit{neighborhood} of $x$), containing $x$, and
require that the proper scoring\break rule~$S(x,\bp)$ be
expressible as a~function of $x$ and the restriction $\bpp_{N_x}$ of~$\bpp$ to
$N_x$:
\[
S(x,\bp)= s(x,\bpp_{N_x}).
\]
We say that such a~scoring rule is \textit{$\cN$-local}, where $\cN=\{
N_x\dvtx x\in\cX\}$ is the neighborhood system. Similarly, its
$0$-homogeneous extension is said to be \textit{$\cN$-local} if
$S(x,\balpha)= s(x,\balpha_{N_x})$. Note this property is strictly
stronger; see Section~\ref{seclogscore} below.

Suppose that the $0$-homogeneous extension of $S$ is continuously
differentiable and $\cN$-local. We then obtain from~(\ref{eqd2}) that,
if $x\notin N_y$,
\[
\frac{\partial S(x, \balpha)}{\partial\alpha_y} =
\frac{\partial S(y, \balpha)}{\partial\alpha_x}=0,
\]
so that without loss of generality we can also require $y \notin
N_x$. Hence, for scoring rules with $\cN$-local $0$-homogeneous
extensions we can assume that the neighborhood relation is symmetric
and so determined by an undirected graph $\graf$ so that $y\in N_x$ if
and only if $x=y$ or $x$--$y$, that is, $x$ and $y$ are neighbors in
$\graf$. We then also say that the scoring rule and its extension are
\textit{$\graf$-local}.

We note that \textit{a~scoring rule with a~$\graf$-local $0$-homogeneous extension only
depends on $\bp$ through its conditional distribution $\bpp_{ |
N_x}$ of $X$ given $X\in N_x$}, that is, it satisfies
%
\begin{equation}
\label{eqcondloc}S(x,P)= s(x,\bpp_{ | N_x}).
\end{equation}
In particular, only knowledge of $\bpp$ up to a~constant factor is
necessary to calculate $S(x,P)$.
Conversely, \textit{the $0$-homogeneous extension of any scoring rule
satisfying}~(\ref{eqcondloc}) \textit{is $\graf$-local.}

\subsection{Logarithmic score}
\label{seclogscore} The simplest case of a~local scoring rule is\break
where~$S(x,\bp)$ is a~function only of $x$ and $p_x$, and is thus
$\graf_0$-local for the totally disconnected graph $\graf_0$. It is
well known [\citet{Bernardo1979ex}] that (for $\#{{\cal
X}} > 2$) a~scoring rule with this property is proper (and is then
strictly proper) if and only if it has the form
%
\begin{equation}
\label{eqlogscore}
S(x, \bp) = a(x) - \lambda\ln p_x
\end{equation}
with $\lambda> 0$. For $a(x)=0$ this is known as the \textit{log-score}.

As described previously, any scoring rule has a~$0$-homogeneous
extension which in this case is $a(x)
-\lambda\ln\alpha_x + \lambda\ln\alpha_+$; however, the extension
depends, not just on
$\alpha_x$, but on $\alpha_y$ for all $y\in{\cal X}$. Hence,
although the scoring rule itself is local, its $0$-homogenous extension
is not, reflected in the fact that knowledge of $\bpp$ up to a~constant factor is not sufficient for calculating the log-score. In
fact, there is no nontrivial proper scoring rule with a~$\graf
_0$-local $0$-homogeneous extension.

Note that the (Shannon) entropy function $H(\balpha)=-\lambda\sum
_x\alpha_x \log(\alpha_x/\alpha_+)$ for the log-score is regular
although the log-score itself is unbounded.

\subsection{Additive scoring rules}
\label{secgraphmod}

Here we describe a~simple way of constructing a~$0$-homogeneous
local scoring rule.
Let $\cB$ be a~collection of subsets of~$\cX$, define ${\cal A}_B
\defeq\{\balpha_B\dvtx
\balpha\in{\cal B}\}$, and let $H_B\dvtx {\cal
A}_B\rightarrow\reals$
be a~concave and $1$-homogeneous function of $\balpha_B$---and thus
also, by extension of its domain, of~$\balpha$. Let $\nabla
H_B\in{\cal A}_B$ be a~$0$-homogeneous supergradient of $H_B$ on
${\cal A}_B$; this is also a~$0$-homogeneous supergradient on the
extended domain ${\cal A}$, if we define its components for $x\notin
B$ as $0$. By the results of Section~\ref{secconstruct}, this determines a~proper scoring rule $S_B(x,\balpha)$. Moreover, $S_B$ vanishes if
$x\notin B$, and otherwise depends on $\balpha$ only through
$\balpha_B$.
We now let
%
\begin{equation}\label{eqaddscore} S(x,\balpha) = \sum_{B\in{\cal
B}} S_B(x,\balpha_B),\qquad  H(\balpha) = \sum_{B\in{\cal B}}
H_B(\balpha_B)
\end{equation}
and
these define a~proper and $0$-homogeneous scoring rule and its
associated $1$-homogeneous entropy function. We shall say that a~scoring rule and entropy function satisfying~(\ref{eqaddscore}) are
\textit{$\cal B$-additive}. When each $H_B$ is a~differentiable function
of~$\balpha_B$, the gradient of $H$ will be the unique associated
scoring rule $S$ of form~(\ref{eqaddscore}).

We note that if we
define an undirected graph $\graf$ by $x$--$y$ if and only if $x,y\in
B$ for some $B\in\cB$, we have that the ($0$-homogenous extension of)
\textit{any $\cB$-additive scoring rule is $\graf$-local}.
If~${\cal C}$ denotes the collection of all cliques of~$\graf$, that
is, all maximal complete subsets of $\cX$, we can collect terms
appropriately and rewrite the expansions in~(\ref{eqaddscore}) above as
%
\begin{equation}
\label{eqgrafscore}
S(x,\balpha) = \sum_{C\in{\cal C}} s_C(x,\balpha_C)
\end{equation}
and
%
\begin{equation}
\label{eqhsum}
H(\balpha) = \sum_{C\in{\cal C}} h_C(\balpha_C).
\end{equation}
We say that a~scoring rule $S$ and entropy function $H$ having the
forms of~(\ref{eqgrafscore}) and~(\ref{eqhsum}) are
\textit{$\graf$-additive.}

We remark that the above constructions can also be applied
straightforwardly to the case of a~countably infinite sample space
${\cal X}$, so long as every set $B\in\cB$ is finite.

We shall in Section~\ref{secconverse} give conditions for the converse to
hold, that is, conditions for a~$\graf$-local scoring rule to be
$\graf$-additive as above, without necessarily demanding each term of
the decomposition~(\ref{eqhsum}) to be concave or $1$-homogeneous.

\section{Examples}
\label{secex}
This section gives some examples of $\graf$-additive and $\graf
$-local scoring rules.

\subsection{Local scoring rules for integer-valued outcomes} We first
consider cases where the outcomes are nonnegative integers.
%
\begin{expl}[(Pair scoring rule)]
Suppose $\mathcal{X}=\{0,1,2,\ldots\}$, and let the graph $\graf$
have edges between successive integers.
The cliques are just the pairs, $C_x \defeq\{x, x+1\}$ $(x=0, 1,
\ldots)$, and a~concave, $1$-homogeneous local entropy function on
$C_x$ has the form $H_x(\alpha_x, \alpha_{x+1}) = \alpha_x
G_x(\alpha_{x+1}/\alpha_x)$ with~$G_x$ concave. The associated
additive scoring rule is
\[
S(x, \bp) = G_{x-1}'\biggl(\frac{p_x}{p_{x-1}}\biggr) +
G_x\biggl(\frac{p_{x+1}}{p_x}\biggr) -
\frac{p_{x+1}}{p_x}  G_x'\biggl(\frac{p_{x+1}}{p_x}\biggr)\qquad
(x=0,1,\ldots)
\]
with the first term absent if $x=0$. The total score based on a~sample $(x_1, \ldots, x_n)$ in which the frequency of $y$ is $f_y$
$(y=0,1,\ldots)$ is thus
\[
\sum_{y=0}^\infty f_y G_y(v_y) + (f_{y+1}-f_y v_y) G_y'(v_y)
\]
with $v_y \defeq p_{y+1}/p_y$. If, for example, we wished to fit the
Poisson model $p_x \propto\theta^x/x!$, we could estimate $\theta$ by
minimizing the total empirical score
\[
\sum_{y=0}^\infty f_y G_y\biggl(\frac{\theta}{y+1}\biggr) +
\biggl(f_{y+1}- \frac{f_y}{y+1} \theta\biggr) G_y'\biggl(\frac
{\theta}{y+1}\biggr).
\]
Taking $G_x(v) = - (x+1)^a~v^{m}/m(m-1)$ for $m\neq0,1$, we obtain the
unbiased estimating equation
%
\begin{equation}
\label{eqesteq}
\theta\sum_{y\geq0} \frac{f_y}{(y+1)^{m-a}} - \sum_{y\geq0}
\frac{f_{y+1}}{(y+1)^{m-a-1}}= 0
\end{equation}
yielding a~simple explicit formula for the estimate. When $m=a$, we
recover the maximum likelihood estimate
$\hat{\theta}=\bar{x}$.
\end{expl}
%
\begin{expl}[(Capture--recapture)]
Consider the following
experiment performed to estimate the number $N$ of fish in a~lake. On $c$
consecutive occasions we catch a~fish, at random, and then replace it.
When a~fish is first caught it is given a~unique tag, so that it can
be recognized on recapture. Each fish $i = 1, \ldots, N$ in the lake
has an associated random variable $X_i$, the number of times it is
caught. For large $c$ and $N$ we can approximate the distribution of
$X_i$ by the Poisson distribution with mean $\theta= c/N$. We will
know $f_x$, the number of fish caught $x$ times, for $x > 0$, but not
$f_0$, the number of fish never caught. The observed data thus arise
from a~truncated Poisson distribution for $X$, conditioned on
$X>0$. If we can estimate $\theta$ we can estimate $N = c/\theta$.
However, because of the need to work with the normalization constant
of the truncated Poisson distribution, the maximum likelihood estimate
of $\theta$ cannot be expressed in explicit form and must be
determined numerically.

Homogeneous local scoring rules can be used to avoid the normalization
constant problem and obtain an explicit estimate. We simply modify
the above analysis of the full Poisson model by removing the edge
$0 $--$1$ from the neighborhood graph $\graf$, together with its
associated local entropy function. Equivalently, we redefine $G_0
\equiv0$. With the other explicit choices made above, the resulting
estimating equation is given by~(\ref{eqesteq}) but with the sums now
over $y \geq1$. For $m=a$ we obtain
%
\begin{equation}\label{eqscorecap}
\tilde{\theta}=\frac{c-f_{1}}{n},
\end{equation}
where $n=\sum_{x\geq1}f_{x}$ is the number of different fish caught.
Note that $c-f_1$ is the number of times a~catch yields a~fish which is
already marked.
In comparison the maximum likelihood estimate $\hat\theta$ satisfies
\[
\hat\theta=\frac{c-ne^{-\hat\theta}}{n},
\]
so $f_1$ in~(\ref{eqscorecap}) is replaced by the estimate of its
expectation $\E_{\hat\theta}(F_1)=ne^{-\hat\theta}$.

In the interests of robustness we might also omit other data, and
again this is easily done. For example, let the only edge in $\graf$
be $1 $--$2$; equivalently, we take $G_x = 0$ for $x \neq1$. Then,
as well as $f_{0}$, the counts $f_{3},f_{4},\ldots$ are also excluded,
only the terms for $y=1$ remain in~(\ref{eqesteq}), and we obtain the
robust Zelterman estimate [\citet{zelterman88}]:
\[
\check{\theta}=\frac{2f_{2}}{f_{1}}.
\]
\end{expl}

\subsection{Local scoring rules for product spaces}
\label{secmulti}

Suppose our discrete sample space is itself a~product space, ${\cal X}
= {\cal X}_1 \times{\cal X}_2 \times\cdots\times{\cal X}_k$. A~point $x$
of ${\cal X}$ has the form $x = (x_1, \ldots, x_k)$. We can define a~useful symmetric neighborhood relation on ${\cal X}$ by $x $--$y$
if, for some $i$, $x^{\setminus i} = y^{\setminus i}$, where
$x^{\setminus i} \defeq(x_j\dvtx j \neq i)$. A~maximal clique of the
associated graph $\graf$ is then defined by an index $i\in\{1,
\ldots,k\}$ and a~vector $\xi^{\setminus i}\in{\cal X}^{\setminus i}
\defeq\bigtimes_{j \neq i} {\cal X}_j$, and has the form
$C_{i,\xi^{\setminus i}} := \{x\dvtx x^{\setminus i} = \xi^{\setminus
i}\}$. Within such a~clique, only the value of $x_i$ can vary, over
the space ${\cal X}_i$.

We can introduce, for such a~clique $C = C_{i,\xi^{\setminus i}}$, a~$1$-homogeneous concave function $H_C$ of $\balpha_C$. Its gradient
will determine a~$0$-homogeneous proper scoring rule $S_C(x,
\balpha)$, vanishing unless $x \in C$, that is, $x^{\setminus i} =
\xi^{\setminus i}$, and in this case depending on~$\balpha$ only
through~$\balpha_C$. In particular, $S_C(x, \bp)$ depends on~$\bp$
only through the implied conditional distribution $\bp(  \cdot\mid
X^{\setminus i} = \xi^{\setminus i})$, for~$X_i$, given $X^{\setminus
i} = \xi^{\setminus i}$.

Conversely, any proper scoring rule
defined for outcomes in ${\cal X}_i$ and distributions over ${\cal
X}_i$ can be applied to the observed value $x_i$ of $X_i$ and the
conditional distribution $\bp(  \cdot\mid X^{\setminus i} =
\xi^{\setminus i})$ (and taken as $0$ if $x^{\setminus i} \neq
\xi^{\setminus i}$): when denormalized, this will be of the above
form. The general $\graf$-additive scoring rule can then be formed
by aggregating a~collection of such single-clique component scoring
rules:
%
\begin{equation}
\label{eqgengraf}
S(x, \bp) = \sum_C S_C\{x_i, \bp(  \cdot\mid
X^{\setminus i} = \xi^{\setminus i})\}  1(x^{\setminus i} =
\xi^{\setminus i}).
\end{equation}

\subsubsection{Specialization}
\label{secspec}

Although it is allowable that the form of the component scoring rule
$S_C$ in~(\ref{eqgengraf}) might vary with the conditioning values
$\xi^{\setminus i}$ that, together with the index $i$, determine the
clique $C$, this level of generality will rarely be needed, and we
might thus restrict attention to proper scoring rules of the form
%
\begin{equation}
\label{eqnongen}
S(x, \bp) = \sum_{i=1}^k S_i\{x_i, \bp(  \cdot\mid
X^{\setminus i}=x^{\setminus i})
\},
\end{equation}
where $S_i$ is a~proper scoring rule for variables in, and
distributions over, ${\cal X}_i$. The associated discrepancy function
is
\[
d(\bp, \bq) = \sum_i \E_{X \sim\bp}  d_i\{ \bp(  \cdot
\mid
X^{\setminus i}), \bq(  \cdot\mid X^{\setminus i})\},
\]
where $d_i$ is the discrepancy function associated with $S_i$.

Recall that we are assuming that $\bp, \bq$ are everywhere positive
distributions. If now each $S_i$ is strictly proper, then $d(\bp,
\bq) = 0$ if and only if, for all $i$ and $x^{\setminus i}$,
$\bq(  \cdot\mid X^{\setminus i} = x^{\setminus i}) = \bp(
\cdot
\mid X^{\setminus i} = x^{\setminus i})$. But (with strict
positivity) this can only occur if $\bq= \bp$, so $S$ is strictly
proper.

\subsubsection{Markov models}
\label{secmarkmod}

Suppose now that we have an undirected graph ${\cal K}$ with vertices
$\{1, \ldots, k\}$, and we restrict attention to distributions
$\bp$ that are Markov with respect to ${\cal K}$. Any component score
$S_C$ in~(\ref{eqgengraf}), or $S_i$ in~(\ref{eqnongen}), will then depend
only on the value $x_i$ of $X_i$, and the conditional distribution of
$X_i$ \textit{given the neighbors} $X^{\neighbours{i}}$ of $i$ in
${\cal K}$. It is possible to calculate this conditional distribution
without having access to the normalizing constant of the overall
distribution $\bp$, which is often hard to compute. In particular, in
the estimation context described in the
\hyperref[secintro]{Introduction}, this can greatly
ease construction and solution of the unbiased estimating equation
associated with this scoring rule.

A~prominent example of a~scoring rule of this kind is the
\textit{pseudo-likelihood} function introduced by \citet{besag75}:
%
\begin{expl}[(Pseudo-likelihood)]
When every component scoring\break rule~$S_i$ in~(\ref{eqnongen}) is the log
score, the overall rule will be just the negative logarithm, $S(x,
\bp) = -\log\mbox{PL}(\bp,x)$, of the pseudo-likelihood
function, defined, for a~joint distribution $\bp$ and
outcome vector $x$, as
\[
\mbox{PL}(\bp,x) \defeq
\prod_{i} \bp(X_{i}=x_{i} \mid X^{\setminus
i}=x^{\setminus i})
\]
(where
in the context of a~Markov model the conditioning variables
$X^{\setminus
i}$ can be reduced to $X^{\neighbours{i}}$). Hence general properties
of proper scoring rules can be applied to pseudo-likelihood. In
particular, a~maximum pseudo-likelihood estimator will typically
be consistent under independent and identically distributed repetitions
(though this argument does
not address consistency under increasing dimension $k$, which is
more relevant in many applications of pseudo-likelihood).
\end{expl}

Replacing the log score with the Brier score leads to the method of
\textit{ratio matching} [\citet{Hyvarinenext}].
%
\begin{expl}[(Ratio matching)]
For the case ${\cal X}_i = \{0,1\}$, take
every component score $S_i$ in~(\ref{eqnongen}) to be the Brier
score, leading to the overall scoring rule
%
\begin{equation}
\label{eqbrierscore}
S(x, \bp) = \sum_i \{x_i - \bp(X_{i}= 1 \mid X^{\setminus
i}=x^{\setminus i})\}^2.
\end{equation}
For a~parametric model ${\cal Q} = \{\bq_\theta\dvtx \theta\in
\Theta\}$, we could estimate $\theta$ by minimizing $\sum_i S(x_i,
\bq_\theta)$. This would equivalently minimize the empirical
discrepancy $d(\hat\bp_n, \bq_\theta)$, where
\[
d(\bp,\bq) =
\sum_{i=1}^k  \sum_{\xi^{\setminus i}\in\{0,1\}^{k-1}}
\bp(X^{\setminus i}=\xi^{\setminus i}) \{
\bp^{
\xi^{\setminus i}}(X_i = 1)- \bq^{
\xi^{\setminus i}}(X_i = 1)
\}^2 \label{eqbrierdisc}
\]
with $\bp^{
\xi^{\setminus i}}(X_i = 1)=\bp(X_i = 1 | X^{\setminus i} =
\xi^{\setminus i})$, etc.
This can be shown to agree with the more complex formula (13) of
\citet{Hyvarinenext}.\setcounter{footnote}{1}\footnote{The further analysis in that paper
does not agree with our~(\ref{eqbrierscore}), and appears to contain
some errors.}
\end{expl}

\section{Characterizing local scoring rules}
\label{secconverse}

Any positive linear combination of the log-score $ -\lambda\ln p_x$
and a~$\graf$-additive score of form~(\ref{eqgrafscore}) will be
$\graf$-local. We now develop a~converse to this result, assuming
henceforth that $S(x, \bp)$ is continuously differentiable on ${\cal
P}$. Under additional conditions on the neighborhood relation $\cN$,
we show that any proper such local scoring rule must be $\graf$-local
for a~suitably defined graph $\graf$ and equal to a~positive linear
combination of the log-score and a~$0$-homogeneous $\graf$-additive score.

We say $x$ is \textit{related} to $y$ and write $x\sim y$, if $x,y\in N_z$
for some neighborhood $N_z\in\cN$. Let $\rho(x)\defeq\{y\dvtx y\sim x\}
$ denote the set of \textit{relatives} of
$x$. Consider now the following condition on the neighborhood system
$\cN$:
%
\begin{cond}
\label{conddisjoint}
There exist $y_1, y_2 \in{\cal X}$ such that, with $\rho_i \defeq
\rho(y_i)$:
%
\begin{eqnarray}
\label{eqsig1}
\rho_1 \cap\rho_2 &=& \varnothing,\\
\label{eqsig2}
\rho_1 \cup\rho_2 &\neq& {\cal X}.
\end{eqnarray}
\end{cond}

Note that in the special case of the trivial neighborhood system $\cN
_0$, that is, $N_x=\{x\}$ for all $x$, this condition is equivalent to
the condition $\#{{\cal
X}} > 2$, as required for the log score to be the only proper $\cN
_0$-local scoring rule; see Section~\ref{seclogscore}.

Assume now that $S$ is a~proper scoring rule. For fixed
$\bp\in\cP$, $S(\bp,\bq)$ is then minimized in $\bq$, subject to
$\bq\in\cP$, at $\bq=\bp$.
Introducing, for each $P$, a~Lagrange multiplier $\lambda(\bp)$ for the
constraint $\sum_x q_x = 1$, we must thus have
%
\begin{equation}
\label{eqnecessary0}
\sum_{x}p_{x} \frac{\partial}{\partial p_{y}}S(x, \bp) + \lambda
(\bp) = 0 \qquad\mbox{for all $y\in\cX$}.
\end{equation}

In the case of a~$0$-homogeneous proper local scoring rule,
we could without loss of generality assume that the neighborhood system
$\cN=\{N_x,x\in\cX\}$ was determined by an undirected graph $\graf
$. In general this is not necessarily the case, as the following
example shows.
%
\begin{expl}
\label{excounter}
A~simple example that does not satisfy the condition is the
neighborhood system determined by the undirected graph $1 $--$2$
$3$, where
$\rho(1) = \rho(2) = \{1,2\}$, $\rho(3) = \{3\}$.
For this graph we can define a~scoring rule as follows:
\begin{eqnarray*}
S(1, \bp) &=& S(2, \bp) = (1-p_1 - p_2)^2,\\
S(3, \bp) &=& (1-p_3)^2.
\end{eqnarray*}
Then $S$ is $\graf$-local, and can easily be shown to be proper (it is
an affine transformation of the
Brier score for the
event $X = 3$). However, its $0$-homogeneous extension is $S(1,
\balpha) = S(2, \balpha) = (\alpha_3/\alpha_+)^2$, $S(3, \balpha) =
\{(\alpha_1 + \alpha_2)/\alpha_+\}^2$, where $\alpha_+ \defeq
\alpha_1
+ \alpha_2 + \alpha_3$. Thus the $0$-homogeneous extension of~$S$ is
not $\graf$-local, and
in particular not $\graf$-additive.
\end{expl}

For neighborhood systems $\cN$ which satisfy Condition~\ref{conddisjoint} we
have the following lemma.
%
\begin{lemma}
\label{lemdisjoint}
Suppose $S$ is proper and $\cN$-local. If Condition~\ref{conddisjoint} holds,
then $\lambda(\bp)$ satisfying~(\ref{eqnecessary0}) is constant on
${\cal P}$.
\end{lemma}
\begin{pf}
For $\cN$-local $S$, condition~(\ref{eqnecessary0}) gives for any
$y\in
{\cal X}$:
%
\begin{equation}
\label{eqmiss}
-\lambda(\bp) = \sum_{\{x\dvtx y\in N_x\}} p_x \,\frac{\partial
S(x, \bp)}{\partial p_y}
\end{equation}
as ${\partial S(x, \bp)}/{\partial p_y}=0$ unless $y\in N_x$.
For any term in the sum in~(\ref{eqmiss}), $S(x, \bp)$, and thus
${\partial S(x, \bp)}/{\partial p_y}$, depends only on $\bpp_{N_x}$,
hence, since $y\in N_x$, only on $\bpp_{\rho(y)}=\{p_z\dvtx z\in\rho
(y)\}$.
Taking $y = y_1$, this implies that
$\lambda(\bp)$ depends only on $\bpp_1 \defeq\{p_z\dvtx z \in
\rho_1\}$; similarly, $\lambda(\bp)$ depends only on $\bpp_2 \defeq
\{p_z\dvtx z \in\rho_2\}$.

By~(\ref{eqsig2}) we can take $w\in{\cal X}\setminus(\rho_1
\cup\rho_2)$. Starting at $\bpp$, consider a~change~$\delta\bpp_1$
to~$\bpp_1$, such that $p_x + \delta p_x \in(0,1)$
($x \in\rho_1)$, and $\delta p^+_1 \defeq\sum_{x\in\rho_1}
\delta p_x \in(p_w-1,\allowbreak p_w)$. Extend the variation $\delta\bpp_1$
to the whole of $\bpp$ by $\delta p_w = - \delta p^+_1$, $\delta
p_x = 0$ (all other~$x$). Then $\bpp+ \delta\bpp\in{\cal P}$.
Since $\lambda(\bpp)$ depends only on $\bpp_2$, which has not
changed, $\lambda(\bpp+ \delta\bpp) = \lambda(\bpp)$. But we can
also express $\lambda(\bp)$ as $\lambda^*(\bpp_1)$, whence~$\lambda^*$ must be constant in an open neighborhood of $\bpp_1$.
It follows that~$\lambda(\bp)$ is constant on ${\cal P}$.
\end{pf}

We now have that, under Condition~\ref{conddisjoint}, any $\cN
$-local proper
scoring rule must satisfy:

For all $\bp\in{\cal P}$ and all $y\in{\cal X}$,
%
\begin{equation}
\label{eqnecessarylam}
\sum_{x}p_{x} \,\frac{\partial}{\partial p_{y}}S(x, \bp) = -\lambda
\end{equation}
for some scalar $\lambda\in\reals$.\vadjust{\goodbreak}

We note that a~particular $\cN$-local solution of
(\ref{eqnecessarylam}) is given by the log-score, $S(x,\bp) =
-\lambda
\ln p_x$. Because~(\ref{eqnecessarylam}) is linear, the general
solution is thus $S = -\lambda\ln p_x
+ S_0$, where, for all $\bp\in{\cal P}$ and all $y\in{\cal X}$,
$S_0$ satisfies the \textit{key equation}:
%
\begin{equation}
\label{eqkey}
\sum_{x}p_{x}  {\partial}/{\partial p_{y}}S(x, \bp) = 0.
\end{equation}
We thus can, and henceforth shall, restrict attention to such
\textit{key local} scoring rules.
We next show that the $0$-homogeneous extension of a~key local scoring
rule is $\graf$-additive for a~suitable undirected graph $\graf$.

Let $H(\bp) \defeq\sum_x p_x S(x, \bp)$ be the associated entropy
function. Then~(\ref{eqkey}) implies $S(y, \bp) = \partial
H(\bp)/\partial p_y$. It follows that
%
\begin{equation}
\label{eqsymm}
\frac{\partial S(x, \bp)}{\partial p_y} =
\frac{\partial S(y, \bp)}{\partial p_x}.
\end{equation}
Hence if $y\in N_x$ but $x \notin N_y$,
${\partial S(x, \bp)}/{\partial p_y}$ must nevertheless vanish. Let
$\graf$ be the undirected graph in which $x$ and $y$ are neighbors
if \textit{both} $x \in N_y$ and $y \in N_x$.
We call $\graf$ the \textit{symmetric core} of $\cN$. Then any key
$\cN$-local proper scoring rule must in fact be $\graf$-local. So we
henceforth confine attention to $\graf$-locality for an undirected
graph $\graf$ and assume that the neighborhoods are determined by
$\graf$ as $N_x=\{x\}\cup\bound(x)$.
%
\begin{lemma}
\label{lemstar}
Under Condition~\ref{conddisjoint}, if $S$ is a~key $\graf$-local
scoring rule,
its $0$-homogeneous extension is
$\graf$-local.
\end{lemma}
\begin{pf}
With $\bp= \balpha/\alpha_+$ we obtain by differentiation that, for
any $y$,
\[
\frac{\partial S(x, \balpha)}{\partial\alpha_y} =
\frac1 {\alpha_+} \biggl\{\frac{\partial S(x, \bp)}{\partial p_y} -
\sum_{x}p_{x} \,\frac{\partial}{\partial p_{x}}S(y, \bp)\biggr\}=
\frac1 {\alpha_+}\,\frac{\partial S(x, \bp)}{\partial p_y}
\]
as~(\ref{eqkey}) and~(\ref{eqsymm}) imply that the second term within
braces vanishes. Hence 
the result follows.
\end{pf}

Before we proceed to show that under Condition~\ref{conddisjoint},
any key
$\graf$-local scoring rule is $\graf$-additive, the following result
is useful.
%
\begin{lemma}
\label{lemsep}
Suppose $x \!\!\not \,$--$y$, Condition~\ref{conddisjoint} holds, and
$S$ is key
$\graf$-local. Then its entropy function $H$ satisfies
\[
\frac{\partial^2 H(\balpha)}{\partial\alpha_x \,\partial\alpha_y}
= 0.
\]
\end{lemma}
\begin{pf}
In this case, by Theorem~\ref{thmprophom}, ${\partial^2
H(\balpha)}/{\partial\alpha_x \,\partial\alpha_y} =
\partial{S(x,\balpha)}/\allowbreak{\partial\alpha_y} = 0$ by Lemma~\ref{lemstar}.
\end{pf}

The following lemma is straightforward:
%
\begin{lemma}
\label{lemsep2}
Let $ f(\balpha)$ be a~twice continuously differentiable function.
Then ${\partial^2
f(\balpha)}/{\partial\alpha_x \,\partial\alpha_y} =0$ if and only if
%
\begin{equation}
\label{eqhh}
f(\balpha) = f(\balpha_W, \alpha_x, \alpha_y^*) + f(\balpha_W,
\alpha_x^*, \alpha_y) - f(\balpha_W, \alpha_x^*, \alpha_y^*)
\end{equation}
for some (and then any) values $\alpha_x^*$, $\alpha_y^*$, where $W
\defeq{\cal X}\setminus\{x,y\}$.
\end{lemma}

We now proceed to establish $\graf$-additivity for any key
$\graf$-local scoring rule:

\begin{theorem}
\label{thmhc}
Let ${\cal C}$ be the set of maximal cliques of the graph $\graf$
satisfying Condition~\ref{conddisjoint}, and let $S$ be a~key $\graf$-local
scoring rule.
Then we can express the scoring rule and associated entropy function as
%
\begin{equation}
\label{eqhc}
S(\balpha) = \sum_{C\in{\cal C}} s_C(\balpha_C), \qquad  H(\balpha)
= \sum_{C\in{\cal C}} h_C(\balpha_C).
\end{equation}
Further, if $H$ is regular, each term $h_C$ in the expansion can be
chosen to be $1$-homogeneous.
\end{theorem}
\begin{pf}
The proof parallels that of the Hammersley--Clifford theorem as given
in \citet{grimmett73}; see also \citet{lauritzen96}, page 36. For a~fixed~$\balpha^*$ and all subsets $A\subseteq\cX$ we define
%
\begin{equation}\label{eqfirstmob}\eta_A(\balpha_A)=H(\balpha
_A,\balpha^*_{\cX\setminus A})
\end{equation}
and note that then $H(\balpha)=\eta_\cX(\balpha)$.
Next, we define for all $B\subseteq\cX$
%
\begin{equation}\label{eqmobdef}h_B(\balpha_B)=\sum_{A\dvtx A\subseteq
B} (-1)^{|B\setminus A|}\eta_A(\balpha_A).
\end{equation}
The M\"{o}bius inversion formula [see, e.g.,
\citet{lauritzen96}, page 239] then yields: for all $A\subseteq\cX$,
\[
\eta_A(\balpha_A)=\sum_{B\dvtx B\subseteq A} h_B(\balpha_B).
\]
Thus, taking $A=\cX$, we have established
%
\begin{equation}\label{eqadditivity}
H(\balpha)=\sum_{B\dvtx B\subseteq\cX} h_B(\balpha_B).
\end{equation}
We next show that all terms $h_B$ in~(\ref{eqadditivity}) vanish unless
$B$ is a~complete set and~(\ref{eqhc}) then follows by collecting
appropriate terms. So suppose there exist $x, y\in B$ with
$x \!\!\not \,$--$y$. We then let $D=B\setminus\{x,y\}$ and write the expression in~(\ref{eqmobdef}) as
\[
h_B(\balpha_B)=\sum_{A\dvtx A\subseteq D} (-1)^{|D\setminus A|}\bigl(\eta
_A-\eta_{A\cup\{x\}}
-\eta_{A\cup\{y\}}+\eta_{A\cup\{x,y\}}\bigr),
\]
where we have abbreviated $\eta_U\defeq\eta_U(\balpha_U)=H(\balpha
_U,\alpha^*_{\cX\setminus U})$, etc. But each of the terms in this
expansion vanishes by Lemma~\ref{lemsep2} and hence $h_B$ vanishes as
required.

If $H$ is regular, we can choose $\alpha^*_x=0$ for all $x\in\cX$ so
each function $\eta_A$ in~(\ref{eqfirstmob}) and hence each $h_B$ in
(\ref{eqmobdef}) will be $1$-homogeneous.

This establishes the desired expansion of the entropy function. The
expansion for the scoring rule is obtained by forming gradients.
\end{pf}

It does not seem to be true in general that each term in~(\ref{eqhc})
can also be chosen to be concave. If this is indeed the case, $S$ and
$H$ are built up additively from proper scoring rules and entropy
functions defined on cliques.

\section{Summary and discussion}
\label{secdisc}
We have defined a~proper scoring rule $S(x$, $P)=s(x,\bpp)$ for discrete
sample space $\cX$ to be local relative to a~neighborhood system $\cN
=\{N_x\}_{x\in\cX}$ if each $s(x,\bpp)$ only depends on $\bpp$
through its restriction $\bpp_{N_x}$ to $N_x$, and shown how to
construct such scoring rules from additive components. Conversely, we
have shown that under appropriate regularity conditions any proper
local scoring rule has this structure, although the additive components
may not in general each correspond to proper scoring rules.

A~definition of homogeneous local scoring rule for a~general
well-behaved topological outcome space $\cX$ that would unify the
discrete and continuous case would be to say that a~scoring rule is
homogeneous and local if it satisfies
%
\begin{equation}\label{eqgenlocal}S(x,P) = S\{x, P(\cdot\mid N_x)\}
\end{equation}
for every open neighborhood $N_x$ of $x$ and every $x\in\cX$.

Clearly, the homogeneous local scoring rules investigated in \citet
{mfpapdsllplsr} satisfy this requirement. It would be interesting to
obtain a~complete characterization of proper scoring rules satisfying
(\ref{eqgenlocal}).

\section*{Acknowledgment}

We are grateful to Valentina Mameli for comments on an
earlier version of this article.



\printaddresses

\end{document}